\documentclass{amsart}
\usepackage{amsmath}
\usepackage{amssymb}
\usepackage{graphicx}
\usepackage{subfigure}

\newtheorem{theorem}{Theorem}[section] 
\newtheorem{lemma}[theorem]{Lemma}
\newtheorem{proposition}[theorem]{Proposition}
\newtheorem{corollary}[theorem]{Corollary}

\theoremstyle{definition}
\newtheorem{definition}[theorem]{Definition}
\newtheorem{example}[theorem]{Example}

\newtheorem{remark}[theorem]{Remark}
\newtheorem{ack}{Acknowledgements}

\newcommand{\we}{\wedge}
\newcommand{\X}{\ensuremath{\mathcal X}}

\renewcommand{\S}{\ensuremath{\mathcal S}}
\newcommand{\B}{\ensuremath{\mathcal B}}

\newcommand{\C}{\ensuremath{\mathbb C}}
\newcommand{\CP}{\ensuremath{\mathbb C}{\mathbb P}}
\newcommand{\N}{\ensuremath{\mathcal N}}
\renewcommand{\L}{\ensuremath{\mathcal {L}}}
\newcommand{\R}{\ensuremath{\mathcal R}}
\newcommand{\Z}{\ensuremath{\mathbb Z}}
\newcommand{\K}{\ensuremath{\mathbb K}}
\newcommand{\Q}{\ensuremath{\mathbb Q}}

\begin{document}
\title{Multinets, resonance varieties, and pencils of plane curves} 
\author{Michael Falk}
\email{michael.falk@nau.edu}
\address{Michael Falk,  Department of
Mathematics and Statistics\\ Northern Arizona University \\ Flagstaff, AZ 86011-5717}

\author{Sergey Yuzvinsky}
\email{yuz@math.uoregon.edu} 
\address{Sergey Yuzvinsky, Department of Mathematics\\ University of Oregon \\ Eugene, OR 94703}

\dedicatory{Dedicated to Joseph H.M. Steenbrink on the occasion of his sixtieth birthday}


\begin{abstract}We show that a line arrangement in the complex projective plane supports a nontrivial resonance 
variety if and only if it is the underlying arrangement of a ``multinet," a multi-arrangement
with a partition into three or more equinumerous classes which have equal multiplicities at each inter-class intersection point, and satisfy a connectivity condition. 
We also prove that this combinatorial structure is equivalent to the existence of a pencil of plane curves, also satisfying a connectivity condition, whose singular
fibers include at least three products of lines, which comprise the arrangement. We derive numerical conditions which impose restrictions on the number of classes, and the line and point multiplicities that can appear in multinets, and allow us to detect whether the associated pencils yield nonlinear fiberings of the complement.
\end{abstract}
\maketitle

\setlength{\parindent}{0pt}
\setlength{\parskip}{.3cm}

\begin{section}{Introduction}

Resonance varieties of complex hyperplane arrangements are the focus of much of the current research in the field, with connections to topics of algebraic, analytic, and topological interest such as cohomology of local systems, generalized hypergeometric functions, lower central series formulae, and Alexander invariants \cite{ESV,SchSuc2,Suc2}. In addition, the description of these varieties relates to classical geometric and combinatorial constructions \cite{Yuz04,Fa04}. 

This paper continues the theme of connections between resonance varieties in degree one and geometry of the complex projective plane. It is based on a combination
of two ideas. In \cite{LY00}, every essential component of $R^1$ was described
by a generalized Cartan matrix $Q$ constructed from the combinatorics of lines.
 The most symmetric particular case of this construction led
to very special configurations of lines, called nets \cite{Yuz04}. On the other
hand, the first author, 
studying the implications of \cite{LY00}  with respect to $K(\pi,1)$ arrangements
arising from pencils of curves, discovered an empirical method to construct from a pair of resonant weights a multi-arrangement of lines which, in a sense, lies within a pencil of curves, resulting in nets ``with multiplicities" \cite{Fa10,Fa11}.

In this paper we call a line arrangement with this combinatorial structure 
a ``multinet." Our main result is that existence of this structure is
equivalent to the existence of a global component of the resonance variety of degree one, and to the existence of a pencil of curves respecting the arrangement (a ``pencil of Ceva type").
The multiplicities of lines come naturally from the combinatorics of
the arrangement.
In some cases these Ceva pencils yield nonlinear fiberings of the complement of the underlying arrangement, consequently the complement is an aspherical space.

Let us establish some notation and give a brief more technical
description of our results. Let $\L=\{\ell_1,\ldots, \ell_n\}$ denote a central arrangement of (homogeneous linear) hyperplanes in $\C^k.$ Let \K\ be a field, and 
let ${\mathcal E}=\Lambda(e_1,\ldots,e_n)$ the exterior algebra over \K\ generated by degree-one elements $e_i$ corresponding to the hyperplanes $\ell_i.$ Define $\partial: {\mathcal E} \to {\mathcal E}$ of degree -1 by $$\partial (e_{i_1} \cdots  \, e_{i_p}) = \sum_{k=1}^p (-1)^{k-1} e_{i_1} \cdots \, \hat{e}_{i_k}\cdots \, e_{i_p},$$ where as usual \  $\hat{}$ \ denotes omission. Let ${\mathcal I}$ be the ideal of $\mathcal E$ generated by $$\{\partial (e_{i_1} \cdots  e_{i_p}) \ | \ \{\ell_{i_1}, \ldots , \ell_{i_p}\} \ \text{is a minimal dependent subset of} \ \L\}.$$ The {\em Orlik-Solomon algebra} $A=A(\L)$ is the quotient $\mathcal E/\mathcal I,$ a graded algebra isomorphic (for $\K=\C$) to the DeRham cohomology of the complement $M=\C^{\ell} - \bigcup_{i=1}^n H_i,$  \cite{OS80}. The generator $\omega_i$ of $A$ represented by $e_i$ is identified under this isomorphism with the holomorphic logarithmic one-form $d(\log(\alpha_i))$, where $\alpha_i$ is a defining linear form for $H_i.$

For the purposes of this paper it will be no loss to assume \L\ is an arrangement in $\C^3.$ Projectivizing we think of \L\ as an arrangement in $\CP^2,$ whose elements are lines, and whose codimension-two intersections are points. Our techniques require that \K\ have characteristic zero, and indeed the results are false over fields of positive characteristic. All of our matrices have integer entries, and there is no loss in assuming $\K=\Q.$

The degree one resonance variety $\R^1(\L)$ is by definition the set of elements
$a\in A^1$ such that each $a$ annihilates some element $b\in A^1$ not proportional to
$a$.
In \cite{LY00}, the second author and A. Libgober gave a description of $\R^1(\L)$ using generalized Cartan matrices. This approach is the basis for the present work, so we will summarize their results in some detail, and refer the reader to \cite{LY00} for proofs. 

Suppose \X\ is a collection of points of multiplicity at least 3. Let $\L'\subseteq \L$ be the set of lines passing through points of \X, and form the $|\X|\times |\L'|$ incidence matrix $J$. 
Let $E$ denote the matrix whose every entry is 1, and let $Q(\X)=J^TJ - E.$ 
The matrix $Q(\X)$ is a generalized Cartan matrix, and so the Vinberg classification
 \cite{Ka} applies. One writes $Q(\X)$ as a block direct sum of indecomposable Cartan matrices, each of which is of finite, affine, or indefinite type. 

It turns out that for any irreducible component of $\R^1$ there exists a set $\X$ of multiple points such that the given component is equal to $\ker(Q(\X))\cap \ker(E)$. Moreover such a
$Q(\X)$ will contain at least three affine blocks, and no finite or indefinite blocks. 
So, let $Q(\X)=Q_1\oplus \cdots \oplus Q_k$ with $k\geq 3$ and each $Q_i$ affine. By the Vinberg classification, $\ker(Q_i)$ is one-dimensional, spanned by a positive integer vector $u_i.$ Then $\ker(Q(\X))\cap \ker(E)=\{\sum_{i=1}^k \lambda_i u_i \ | \ \sum_{i=1}^k \lambda_i=0\}.$ 
One can assume without any loss that the sums of the coordinates of $u_i$ are pairwise equal. Then the vectors $u_i-u_1$ form a basis for $\ker(Q(\X))\cap\ker(E)$.
The block sum decomposition of $Q(\X)$ defines a partition $\N=\L_1'\cup \cdots \cup \L_k'$  of $\L',$ with the non-negative vector $u_i$ supported on the block $\L_i'$ of \N, for $1\leq i\leq k.$

This setup forms the starting point for our work. For convenience let us assume $\L'=\L.$ (In this case we say the corresponding component of $\R^1$ is {\em supported on \L.})
We consider the positive integer coordinates of $u_i$ as multiplicities of 
lines in the corresponding block of \N\ and obtain a {\em multinet}.
A multinet consists of a partition $\Pi=\L_1\cup \ldots \cup \L_k$ of \L\ into $k\geq 3$ {\em classes}, a set of points \X, called the {\em base locus}, and a multiplicity function $m: \L \to \Z_{>0}$ such that, counting multiplicities, each class contains the same number of lines, and each point of \X\ is contained in the same number of lines from each class.
There is a further condition (``connectivity") concerning intra-class intersection points corresponding to the indecomposability of the $Q_i.$ This structure is then equivalent to the decomposition of $Q(\X)$ into $k$ indecomposable blocks of affine type. Thus there is a component of the first resonance variety supported on \L\ if and only if \L\ is the underlying arrangement of a multinet. \footnote{A similar result was discovered independently by M.~Marco \cite{Mar05}.}

The two requirements on multiplicities, namely, that classes contain the same number $d$ of lines, and that inter-class intersection points lie in the same number of lines from each class, are necessary conditions for the existence of a pencil of degree $d$ plane curves with three or more (possibly non-reduced) singular fibers whose set-theoretic union is $\bigcup \L.$ The connectivity condition on multinets implies this pencil has ``connected fibers" in the appropriate sense (see Section~\ref{three}). We prove the converse of this: if $\L=\L_1\cup \cdots \cup \L_k$ is a multinet then there is a connected pencil of plane curves among whose singular fibers are the curves given by $\prod_{i=1}^n \alpha_i^{m(\ell_i)},$ where, again, $\alpha_i$ is the homogeneous linear form defining $\ell_i.$ The connectivity condition, for both multinets and pencils, is crucial - the result is false without this restriction. The resonant 1-forms with coefficient vectors $u_i-u_j$ are seen as pullbacks of 1-forms on $\CP^1$ with three or more points removed, whose pairwise products are {\em a fortiori} zero for dimensional reasons. A {\em net} is a multinet in which all lines have multiplicity one and every point of \X\ is contained in precisely one line from each class. In this case the connectivity condition is vacuous.

 For line arrangements  \L\ with a nontrivial first resonance variety, the existence of an associated pencil of curves, whose singular fibers include the lines of \L\ among their components, was proven in \cite{LY00}. Our result is sharper, identifying the pencil precisely, 
 and expressing the multiplicities of the linear factors of the split fibers  in terms of combinatorics of multinets. This pencil may be viewed as a realization, by the restriction of an algebraic map of projective varieties, of the map of the complement $M$ to a curve  arising from a local system with nonvanishing first cohomology, whose existence is ensured by Arapura's Theorem \cite{Arap97}.\footnote{After receiving a preprint of this paper, A.~Dimca has given an argument based on Arapura's theorem of a more general, but less explicit form of the equivalence of (i) and (iii) of Corollary~\ref{bigcor} \cite{Dimca06}.}

There is enough geometric data encoded in the multinet structure for us to write down a Riemann-Hurwitz type formula for the sum of the euler numbers of the special fibers, in terms of the degree $d$ and multiplicities of points inside and outside the base locus. This imposes some restriction on  the various parameters; for instance if all line multiplicities equal one then a multinet can have at most 5 classes, for any $d.$ Comparing with the euler numbers of the special fibers $\bigcup \L_i,$ this formula also detects the existence of additional singular fibers that are not part of $\bigcup \L.$ When there are no other singular fibers, we say the multinet is {\em complete}. In this case the pencil defines a fibering of the complement $M$ with aspherical base and fiber, showing that $M$ is aspherical, i.e., \L\ is a $K(\pi,1)$ arrangement. In this case,  it is sometimes possible to add lines to \L\ preserving the fibering property. These larger fibered $K(\pi,1)$ arrangements are not themselves multinets. The fibered  $K(\pi,1)$ arrangements $J_d, d\geq 2$ constructed in \cite{FaRa86} arise from a family of nets in this way.  Fiber-type arrangements of rank three also arise from this construction, starting from the unique rank-two net. Using other complete multinets, we are able to produce new examples of fibered arrangements by this method.

Here is an outline of the paper. In Section 2, we define multinets and establish a one-to-one correspondence between components of $\R^1(\L)$ with support equal to \L\ (global components) and multinets on \L. 
In Section 3 we show that multinets correspond to connected pencils of curves with three or more singular fibers that are (non-reduced) products of lines.  In Section 4 we write down the Riemann-Hurwitz type formula for our situation, and discuss the implications for the $K(\pi,1)$ problem. Many examples appear throughout the paper.
\end{section}

\begin{section}{Multinets and resonance varieties}

Let $\L$ 
be an arrangement of lines in the complex projective plane. A point contained in two or more lines of \L\ will be called an {\em intersection point.}  
A point contained in three or more lines will be called a {\em multiple point}.

Let $m: \L \to \Z_{>0}$ be a function which assigns to each $\ell\in \L$  a positive integer multiplicity $m(\ell)$. The pair $(\L,m)$ is  called a {\em multi-arrangement}.

A multi-arrangement is a realization of a rank-three matroid whose ground set contains $m(\ell_i)$ copies of $\ell$ for each $\ell_i \in \L.$ Intersection points correspond to the rank-two flats in the matroid. The Orlik-Solomon algebra of \L, defined in the introduction, depends only on the underlying matroid, and the theory developed in this section generalizes to the more general setting of (loopless) matroids.

\begin{definition}
\label{multinets}
A {\em weak $(k,d)$-multinet} on a multi-arrangement $(\L,m)$ is a pair $(\N,\X)$ where \N\ is a partition of \L\ into $k\geq 3$ classes $\L_1, \ldots \L_k,$ and \X\ is a set of multiple points, such that

\begin{enumerate}
\item $\sum_{\ell\in\L_i} m(\ell)=d,$ independent of $i;$
\item For every $\ell\in \L_i$ and $\ell'\in \L_j,$ with $i\not = j$, 
the point $\ell \cap \ell'$ is in $\X;$ 
\item For each $p\in \X$, 
$\sum_{\ell\in \L_i, p\in \ell} m(\ell)$ is constant, independent of $i$.
\end{enumerate}

A {\em multinet} is a weak multinet satisfying the additional property
\begin{enumerate}
\setcounter{enumi}{3}
\item For $1\leq i\leq k,$ for any $\ell,\ell'\in \L_i,$ there is a sequence of $\ell=\ell_0,\ell_1,\ldots,\ell_r=\ell'$ such that $\ell_{j-1}\cap \ell_j\not \in \X$ for $1\leq j \leq r.$
\end{enumerate}
\end{definition}

\begin{remark}\label{prime} Multiplying all $m(\ell)$ of a $(k,d)$-multinet by 
a positive integer $c$ defines a $(k,cd)$ multinet with the same $\L$ and $\X$.
We will always assume that $d$ is chosen minimally in this sense, that is, that the line multiplicities are mutually relatively prime. 
\end{remark}

Examples appear in Figures \ref{multinets-fig} and \ref{steiner}. Each line of multiplicity $m(\ell)\geq 2$ is so labeled; the partition \N\ is indicated with capital letters $A,B,C.$ 

We will call \X\ the {\em base locus} - it need not include all multiple points of \L. Note that, if $(\N,\X)$ is a weak multinet,  then \X\ is determined by \N, namely $\X=\{\ell \cap \ell' \ | \ \ell \in \L_i, \ell'  \in \L_j, i\not = j\}.$ If $(\N,\X)$ is a multinet, then \X\ determines \N\ as well: construct a graph $\Gamma$ with vertex set \L\ and an edge from $\ell$ to $\ell'$ when 
$\ell \cap \ell' \not \in \X.$ Then, by (ii) and (iv),  the $\L_i$ are the components of $\Gamma.$
We will see below (Remark \ref{refinement}) that every weak multinet can be refined to a multinet with the same base locus.

The third condition says that the number of lines from $\L_i$ passing through $p\in \X$, counting multiplicities, is the same for every $i$. This number is denoted $n_p.$ 
If $n_p=1$ for every $p\in \X$ then $(\N,\X)$ is called a {\em net}. In this case condition (iv) follows from (iii), and $m(\ell)=1$ for every $\ell \in \L.$ The converse of the last statement is false - there are multinets with $m(\ell)=1$ for every
$\ell\in\L$ that are not nets (see Figure \ref{steiner}). 
Combinatorially a $k$-net corresponds to a set of $k-2$ mutually orthogonal latin squares; in particular 3-nets correspond to quasi-groups. 
In studies of components of resonance varieties, nets in $\CP^2$
first appeared implicitly
in \cite{LY00} and explicitly with a partial classification in
\cite{Yuz04}. See also \cite{Kaw05}. 

\begin{figure}
\centering
\subfigure[A (3,2)-net]{\includegraphics[width=5cm]{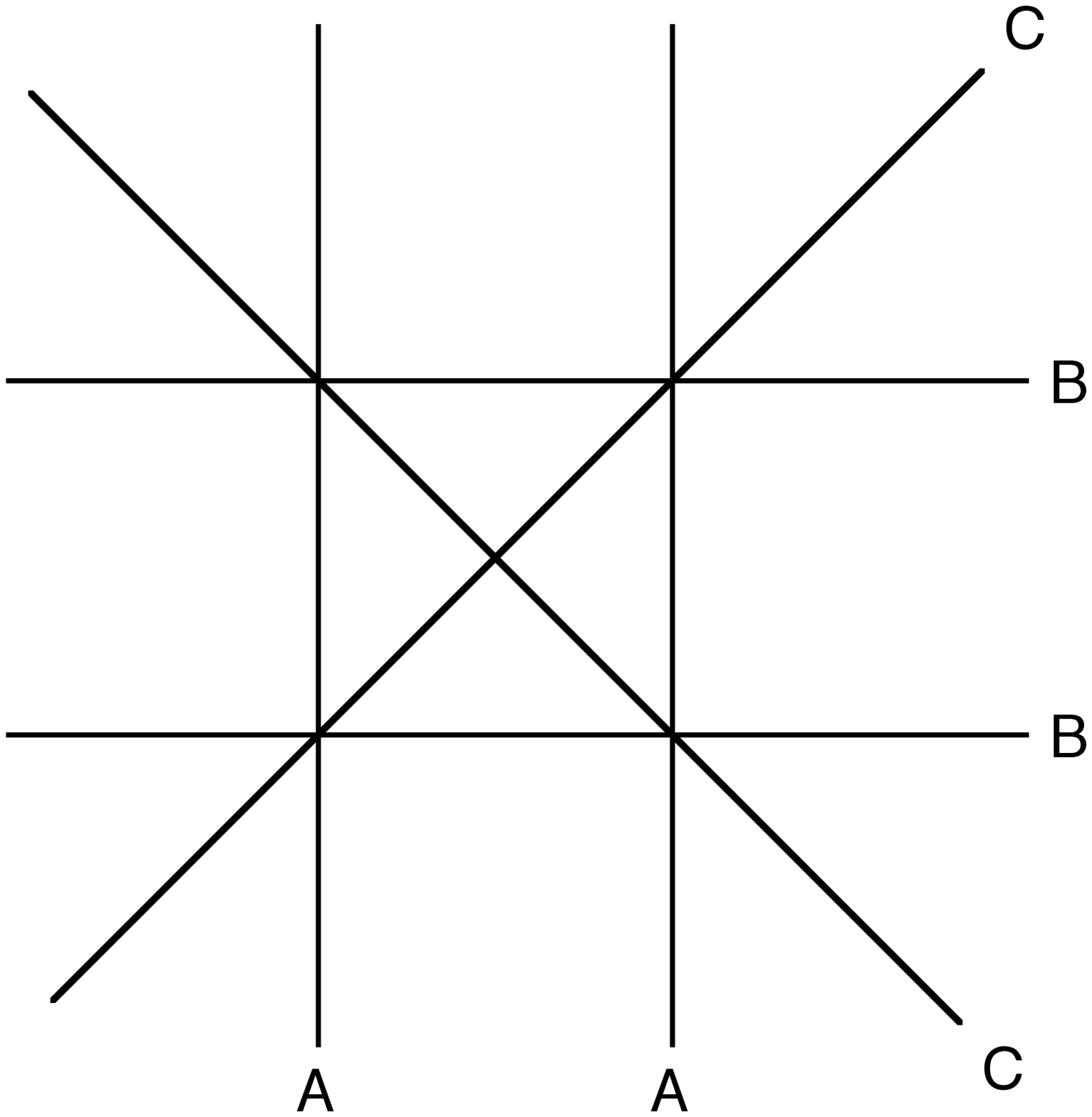}}
\subfigure[A (3,4)-multinet]{\includegraphics[width=5cm]{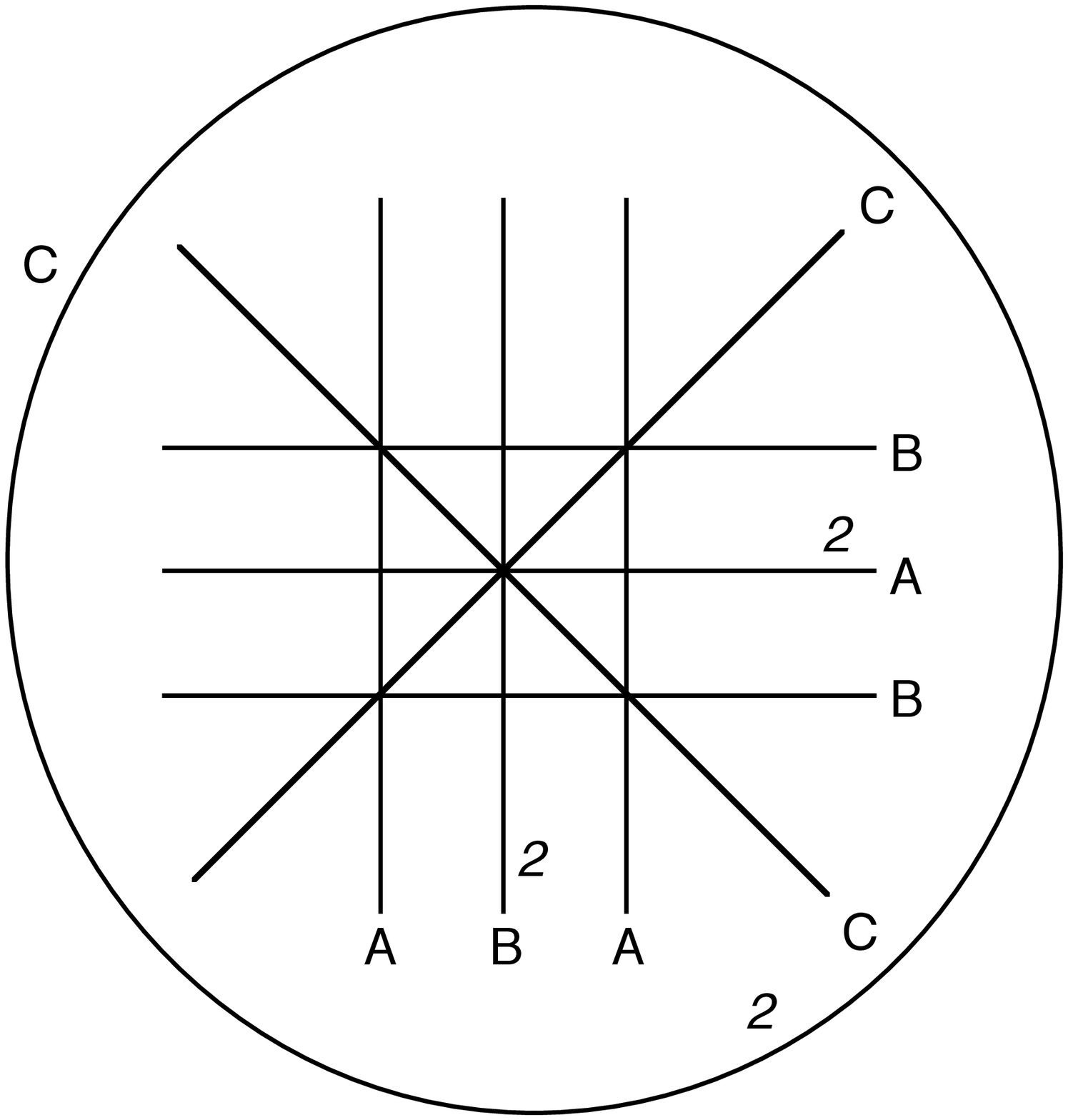}}
\caption{}
\label{multinets-fig}
\end{figure}

\begin{figure} 
\centering
\includegraphics[width=6cm]{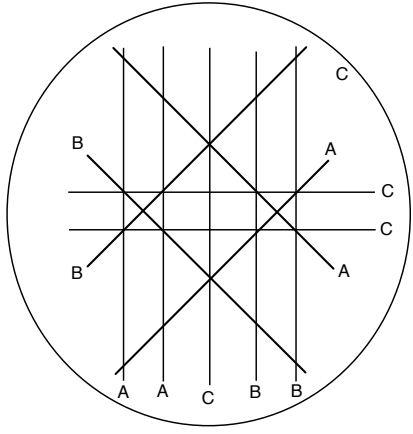} 
\caption{A multinet with $m(\ell)=1$ for all $\ell \in \L,$ which is not a net.}
\label{steiner}
\end{figure}

\begin{lemma} 
\label{numerology}
Suppose $(\N,\X)$ is a weak $(k,d)$-multinet. Then
\begin{enumerate} 
\item $\sum_{\ell\in \L} m(\ell)=dk.$
\item $\sum_{p\in \X} n_p^2=d^2$
\item For each $\ell\in \L,$ $\sum_{p\in\X\cap\ell} n_p=d.$
\end{enumerate}
\end{lemma}

\begin{proof} Each class $\L_i$ consists of $d$ lines counting multiplicity. This implies (i). For (ii), count $|\L_1\times \L_2|,$ with multiplicities, using the fact that every $p\in \X$ contains $n_p$ lines of each of $\L_1$ and $\L_2.$ For (iii), let $\ell \in \L_i,$ and count $\L_j$ with multiplicities, with $j\not = i$. 
Each point $p$ lying on $\ell$ has $n_p$ lines from $\L_j$ 
passing through it, counting multiplicities.
\end{proof}

Given any set \X\ of multiple points, let $J(\X)$ denote the $|\X| \times |\L|$ incidence matrix of \X, whose $(p,\ell)$ entry is 1 precisely when $p\in \ell.$
Let \K\ be a field of characteristic zero, and let $A$ be the Orlik-Solomon algebra of \L\ over \K, as defined in the introduction. 
We denote the generator of $A$ corresponding to $\ell \in \L$ by $\omega_\ell;$ 
then $\{\omega_\ell \ | \ \ell\in\L\}$ is a basis of $A^1$. A subspace $R$ of $A^1$ is {\em isotropic} if $a\we b=0$ for every $a,b \in R.$ Then, perforce, an isotropic subspace of dimension at least two is contained in $\R^1(\L).$

\begin{theorem} Suppose $\L$ supports a weak $(k,d)$-multinet. Then there is a $(k-1)$-dimensional isotropic subspace of $A^1.$
\label{main1}
\end{theorem}

\begin{proof} 
 Let $J=J(\X)$ and $u_i=\sum_{\ell\in \L_i}m(\ell)\omega_\ell.$ Then $Ju_i$ is the vector with entries $n_p, p\in \X,$ for each $1\leq i\leq k,$ by condition (iii) of the definition.
Thus $u_i-u_1$ lies in the kernel of $J$ for each $i.$ 
It follows from section 3 of 
\cite{LY00} that the vectors $u_i-u_1$ have pairwise 
products equal to zero in $A.$ Alternatively, the partition \N\ is neighborly, by (iii), and the vectors $u_i-u_1$ restrict to parallel vectors on each block $\L_i,$ so the pairwise products vanish by the criterion of \cite{Fa97}.
They are linearly independent since the vectors $u_i$
have pairwise disjoint supports.
\end{proof}

Since the criterion of \cite{Fa97} is characteristic-free (see also \cite{Fa04}), the previous theorem also holds for fields of characteristic $p>0$, 
provided each class has a line whose multiplicity is not divisible by $p$ (so that each $u_i\not=0$).

Next we will prove the converse to Theorem \ref{main1}. For this argument it is essential that \K\ have characteristic zero. In this case it was shown in \cite{LY00} that $\R^1(\L)$ is a union of linear subspaces which pairwise intersect trivially. Each of these subspaces is an irreducible
component of $\R^1(\L)$ and is isotropic. In particular the dimension of the subspaces is at least 2.
A component $R$ of $\R^1(\L)$ is called a {\em global resonance component} if $L$ is not contained in any coordinate hyperplane in $A^1.$ We say that \L\ {\em supports} $R$.

\begin{theorem} Suppose \L\ supports a global resonance component $R$ of dimension $k-1.$ 
Then $\L$ supports a $(k,d)$-multinet for some $d.$
\label{main2}
\end{theorem}

\begin{proof} We apply again the results of \cite{LY00}, 
referring the reader to the brief sketch in the introduction.
 The hypothesis implies the existence of a set of multiple points $\X$, 
with incidence matrix $J=J(\X)$, such that the kernel of the matrix $Q=J^TJ-E$
 has dimension $k\geq 3$. Here $E$ is the matrix of all ones. 
The matrix $Q$ is a generalized Cartan matrix,
symmetric with nonnegative integers on the main diagonal and 
-1 or 0 off the main diagonal. There is a
block direct-sum decomposition $Q=Q_1\oplus \cdots \oplus Q_k$ with 
each $Q_i$ indecomposable of affine type. (The absence of finite types 
follows from the hypothesis that $R$ is supported on $\L$.)
The kernel of each $Q_i$ is one-dimensional, spanned by a positive integer
vector $u_i$. Define the class $\L_i$ to be the support of the block $Q_i$. Then $\N=\{\L_1,\ldots,L_k\}$ is a partition of $\L$. We identify $u_i$ with the non-negative integer vector $(u_{i1},\ldots,u_{in}),$ with support $\L_i,$ whose nonzero entries correspond to those of $u_i.$

The integer vectors $u_i$ can be chosen to be primitive vectors.
However, we need to ensure that the sums $\sigma_i=\sum_{j=1}^n u_{ij}$ of their coordinates 
are equal to each other. For that it suffices to multiply each $u_i$ by
$d/\sigma_i$ where $d$ is the least common multiple of all $\sigma_i$.
Then the $u_i$ will be positive integer 
vectors, not necessarily primitive, having the same coordinate sum $d$. Assign multiplicities to
the lines in $\L_i$ to match the corresponding entry in $u_i$, 
for each $i$. Then, by construction, condition (i) of Definition \ref{multinets} holds. The resulting multiplicities $m(\ell), \ell \in \L,$ will be mutually relatively prime.

Suppose $\ell$ and $\ell'$ lie in different classes. Then the $(\ell,\ell')$ 
entry of $Q$ is zero, whence the $(\ell,\ell')$ entry of $J^TJ$ is 1. 
This implies some $p\in \X$ contains both $\ell$ and $\ell',$ i.e. $\ell\cap\ell'\in\X$.

Then, for every $i,j,$ the vector $u_i-u_j$ lies in the kernel of $E$, hence also in the kernel of $J,$ since $\ker(Q)\cap\ker(E)=\ker(J)\cap\ker(E).$ This implies that condition (iii) holds. 
Thus $(\N,\X)$ is a weak $(k,d)$-multinet.

To prove (iv) we study more carefully the formation of classes $\L_i.$ Given the matrix $Q$, form the graph $\Gamma$ with vertex set \L\ and an edge connecting $\ell$ and $\ell'$ if and only if the $(\ell,\ell')$ entry of $Q$ is nonzero. The indecomposable blocks $Q_1,\ldots, Q_k$ of $Q$ are precisely the restrictions of $Q$ to the connected components of $\Gamma.$ 
Now the $(\ell,\ell')$ entry of $Q$ is zero if and only if the same entry of $J^TJ$ is 1, which means that $\ell \cap \ell' \in \X.$ 

So, suppose $\ell, \ell'\in \L_i$ for some $i$. Then the pair $(\ell,\ell')$ indexes an entry in $Q_i,$ an indecomposable component of $Q$.  By the observations in the preceding paragraph, there is a path in $\Gamma$ from $\ell$ to $\ell'$, so there is a sequence $\ell=\ell_0, \ell_1,\ldots, \ell_r=\ell'$ such that $\ell_{j-1} \cap \ell_j \not \in \X$
for $1\leq j\leq q.$ By condition (ii), already proven, each $\ell_j$ is in $\L_i$. Thus (iv) holds, and $(\N,\X)$ is a multinet.
\end{proof}

\begin{remark} Theorems \ref{main1} and \ref{main2} have as
 a consequence that any weak multinet can be refined to a
 multinet with the same base locus $\X$. Indeed, suppose $(\N,\X)$ is a weak $(k,d)$-multinet. Set $J=J(\X)$. Let $w_i=\sum_{\ell \in \L_i} m(\ell)$. Then, for each $i,$  $Jw_i$ has entries $n_p$ for $p\in \X,$ and $J^TJw_i$ is then the vector of all $d$'s, by Lemma \ref{numerology}(iii).  Thus $w_i$ is a non-negative vector in the kernel of $Q=J^TJ-E$, with support $\L_i$, for each $i$. Then $w_i$ is a positive linear combination of the vectors $u_i$ from the proof of Theorem \ref{main2}. It follows that the $(k',d')$-multinet $(\N',\X')$ constructed from $Q$ as in that proof satisfies $d'\leq d, \ \X'=\X,$ and $\N'$ is a refinement of \N.

For example, consider an arrangement of five concurrent lines, one having multiplicity two. Partition the lines into three classes, of two lines each, counting multiplicities. This is a weak $(3,2)$-multinet. Its multinet refinement is the partition into five singletons, with each line having multiplicity one, i.e., a $(5,1)$-net. In this case the matrix $Q$ is the $5\times 5$ zero matrix, the $Q_i$ are $1\times 1$ zero matrices - these are of affine type, with $u_i=1$ for each $i.$ 
\label{refinement}
\end{remark}


\end{section}

\begin{section}{Multinets and pencils of plane curves}
\label{three}
In this section we relate the notion of multinet in $\CP^2$
to pencils of plane algebraic curves. In light of the connection with resonance established in the preceding section, this construction may be viewed as a concrete realization, via the restriction of an algebraic map of projective varieties, of the mapping of the complement of $\L$ to a curve predicted by the theorem of Arapura \cite{Arap97} (generalizing the Castenuovo-de Franchis lemma \cite{GH}), arising from a pair of non-proportional one-forms with wedge product zero. The existence of such a pencil, whose singular locus includes the set-theoretic union of the lines of \L, was proved in \cite{LY00},  and was used in \cite{LY00}  and \cite{Yuz04}  to derive restrictions on nets. Our result is more precise, showing that the multiplicities of lines in the singular locus match the coefficients in the one-forms, and are determined combinatorially. We give a different, more direct proof. In particular our argument does not depend on Arapura's theorem.

We will identify a homogeneous polynomial in three variables (usually determined up to a nonzero multiplicative constant) with the projective 
plane curve it defines, and often refer to either as
 a ``curve." A one-dimensional linear system of curves is called a {\em pencil.}
 One can think of a pencil as a line in the projective space 
of homogeneous polynomials of some fixed degree. Thus any two distinct curves generate a pencil, and conversely a pencil is 
determined by any two of its curves $C_1, C_2$. An arbitrary curve in the 
pencil is then $aC_1 + bC_2,\ [a:b]\in \CP^1.$ Every two curves in a pencil 
intersect in the same set of points $\X=C_1\cap C_2,$ called the base of the pencil.
We will always assume that our pencils have no fixed components, that is,
 that the base locus is a finite set of points.

The two curves $C_1, C_2$ 
determine a rational map $\pi: \CP^2 \rightarrowtail \CP^1$ via $p\mapsto [C_2(p):-C_1(p)]$
whose indeterminacy locus is the base
of the pencil. The (closure of the) fiber of $\pi$ over $[a:b]$ is the curve
$aC_1+bC_2,$ and each point outside the base locus lies in a unique such curve. The map $\pi$ is uniquely determined by the pencil, up to linear change of coordinates in $\CP^1$.
For convenience we will often call $\pi$ a ``pencil," when no confusion will result. For more about pencils the reader may consult \cite{GH}.

A curve of the form  $\prod_{i=1}^q \alpha_i^{m_i},$ where $\alpha_i$ 
is a linear form and $m_i$ is a positive integer, 
for $1\leq i\leq q,$ will be called {\em completely reducible}. We are mainly interested in pencils that have some completely reducible fibers.
For later use we treat the two-dimensional case. We will continue to call the fibers ``curves," thinking of them as unions of lines through the origin in the affine plane. The last assertion below will be used in the next section to detect ``hidden" singular fibers  (Proposition~\ref{localtest}.)

\begin{lemma} Suppose $C_1$ and $C_2$ are (homogeneous) curves of degree $n$ in two variables, with no common components. Let $\pi: \CP^1 \to \CP^1$ be the rational map determined by $C_1$ and $C_2.$ Then 
\begin{enumerate}
\item every fiber of $\pi$ is completely reducible.
\item $\pi$ is a regular map, in particular distinct fibers are disjoint.
\item the generic fiber of $\pi$ consists of $n$ distinct points of multiplicity one.
\item $2n-2=\sum_{p\in \CP^1} (m_p-1),$ where $m_p$ is the multiplicity of $\pi$ at $p.$
\end{enumerate}
\label{rank2ceva}
\end{lemma}

\begin{proof} Any homogeneous polynomial in $\C[x,y]$ is a product of linear forms, so every fiber is completely reducible. Since there are no fixed components, the pencil is base-point free and the associated rational map is regular. It defines a branched covering of $\CP^1$ by $\CP^1$ of degree $n$, so the generic fiber consists of $n$ points of multiplicity one. The Riemann-Hurwitz formula for maps of curves implies $\chi(\CP^1)=n\chi(\CP^1)-\sum_{p\in\CP^1} (m_p-1),$ which is equivalent to (iv).
\end{proof}

We consider pencils in $\CP^2$ generated by two completely reducible curves. We are grateful to M.~Marco for pointing out an error in the first version of Lemma~\ref{nodes}, necessitating the additional hypothesis.

\begin{lemma} Suppose $\pi$ is a pencil with no fixed components and at least two completely reducible fibers $C_1$ and $C_2.$ Let $p\in C_1 \cap C_2$ be a base point, and suppose $C_1$ and $C_2$ have the same multiplicity $n_p$ at $p.$ Then
\begin{enumerate}
\item no two fibers of $\pi$ are tangent at $p$,
\item if $n_p=1$ then the generic fiber of $\pi$ is nonsingular at $p$, and
\item if $n_p>1$ then the generic fiber of $\pi$ has an ordinary singularity of multiplicity $n_p$ at $p$.
\end{enumerate}
\label{nodes}
\end{lemma}

\begin{proof} The tangent cones to $C_1$ and $C_2$ at $p$ generate a pencil of curves in two variables, of degree $n_p,$ whose fibers are the tangent cones to the fibers of $\pi$ at $p$. Because $C_1$ and $C_2$ have no common components, this ``local" pencil satisfies the hypothesis of Lemma~\ref{rank2ceva}. The first assertion then follows from (ii), and the second and third statements follow from (iii) of that lemma.
\end{proof}

Let $\varphi: \S \to \CP^2$ be the blow-up of $\CP^2$ at the points of \X. The absence of fibers tangent at a base point (Lemma~\ref{nodes}(i)) implies that the rational map $\pi: \CP^2 \rightarrowtail \CP^1$ lifts to a regular mapping $\tilde{\pi}: \S \to \CP^1.$ The fibers of $\tilde{\pi}$ are the proper transforms of the fibers of $\pi$ under the blow-up $\varphi.$

We say the pencil $\pi$ is {\em connected} if every fiber of $\tilde{\pi}$ is connected. Equivalently, by Lemma \ref{nodes}(ii) and (iii), $\pi$ is connected if and only if there is no (reducible) fiber of $\pi$ which is a union of finitely many proper subvarieties meeting only in the base locus.


\begin{definition} A {\em pencil of Ceva type} (or ``Ceva pencil") is a connected pencil of plane curves (with no fixed components), in which three or more fibers are completely reducible.
\label{Ceva}
\end{definition}

The condition of Definition~\ref{Ceva} means that, as a line in the projective space $\CP^N$ of degree $d$ curves ($N=\binom{d+2}{2}-1$), a Ceva pencil is a trisecant to the subvariety of completely reducible curves, which can be identified with the Chow variety of 0-cycles of degree $d$ in $\CP^2,$ see \cite[Section 4.2]{GKZ94}. 
The terminology comes from the first example below.

\begin{example} Consider the Fermat pencil $ax^d+by^d+cz^d=0, \ [a:b:c]\in \CP^2, \ a+b+c=0.$ There are three singular values: $[1:-1:0], [0:1:-1],$ and $[-1:0:1].$ The corresponding fibers $x^d-y^d,  \ y^d-z^d,$ and $z^d-x^d,$ are completely reducible. The components of each of the singular fibers meet in a single point outside the base locus, so the pencil is connected. The resulting arrangement of $3d$ lines, called the Ceva arrangement (so named in \cite{BHH87}), supports a $(3,d)$-net; the case $d=2$ is Figure~\ref{multinets-fig}(a). This arrangement appears in \cite{LY00, Fa97,L5} in the context of resonance and characteristic varieties, and  in slightly modified form in \cite{FaRa86}, as an example of a nonlinearly fibered arrangement (see Example~\ref{jd}). 
\label{Ceva-ex}
\end{example}

\begin{example} The other ``classical" example of a pencil of Ceva type is the Hesse pencil of cubics $a(x^3+y^3+z^3) + 3bxyz, [a:b]\in \CP^1$ which share the same nine inflection points. 
In this case there are four completely reducible fibers, each of which is a product of three distinct lines, which meet in pairs outside the base locus, so again the pencil is connected.  
The resulting arrangement of twelve lines in $\CP^2$ is called the Hessian
arrangement \cite{OT92}. It supports a $(4,3)$-net, and in fact is the only known example of a line arrangement supporting a $(4,d)$-net for any $d$ \cite{Yuz04}. 
In this case, again, all other fibers are nonsingular. 
This additional property implies that the complement of the Hessian
arrangement, and of the Ceva arrangements of the preceding example, are aspherical spaces. This was the original motivation 
for studying similarly defined ``sharp pencils" in \cite{Fa10}, also called ``Hessian pencils" in 
\cite{Fa11}, and will be taken up in the next section.
\label{Hesse}
\end{example}

\begin{example} 
The connected pencil $[(x^2-y^2)z^2:(y^2-z^2)x^2]$ has three singular fibers $(x^2-y^2)z^2, (y^2-z^2)x^2,$ and $(z^2-x^2)y^2,$ which are completely reducible but 
not reduced. The resulting arrangement of nine lines is the $B_3$ arrangement (Figure \ref{multinets-fig}(b)). The pencil determines a nonlinear fibering of the complement (see Example~\ref{ex2b3}). This example, discovered by the first author \cite{Fa11}, motivated our work.
\label{exb3}
\end{example}

\begin{proposition} \label{pen-net}
A pencil of Ceva type induces a multinet on 
the line arrangement \L\ consisting of the components of 
its completely reducible fibers.
\end{proposition}

\begin{proof} Suppose $\pi$ is a Ceva pencil of degree $d$ curves with completely reducible fibers $C_1,\ldots,C_k.$ For each $i$ let $\L_i$ be the arrangement of lines defined by the factors of $C_i,$ and let $\L=\L_1\cup \cdots  \cup \L_k$. To each $\ell\in \L$ assign the multiplicity $m(\ell)$ to be the multiplicity of the corresponding linear factor in $C_i.$ Then the $\L_i$ form a partition of \L\ because the pencil has no fixed components. Also, property (i) of Definition~\ref{multinets} is clearly satisfied. Let \X\ be the base of the pencil. Then property (ii) is automatic. Observe that the multiplicity of any fiber $C=aC_1+bC_2$ at $p\in \X$ is at least the minimum of the multiplicities of $C_1$ and $C_2$ at $p.$ Since any two of the $C_i$ generate the pencil, this implies they all must have equal multiplicities at $p.$ This confirms condition (iii).

It is only left to prove that the weak multinet we obtain is really a
multinet. This follows from the connectedness of the pencil. Indeed suppose
that the property (iv) fails. Then the components of the  graph $\Gamma$ determine finitely many proper subvarieties of of the fiber $C_i$ which pairwise intersect only in the base locus. 
Then the pencil $\pi$ is not connected, by our earlier observation.
\end{proof}

Here is an example of a pencil which satisfies the first, but not the second condition of Definition~\ref{Ceva}, and yields a weak multinet that is not a multinet.
\begin{example} Let $C_1=x^d$ and $C_2=y^d.$ Then every curve $aC_1+bC_2$ except
 $C_1$ and $C_2$ is the product of $d$ distinct lines. All these fibers are
completely reduced, with components meeting only inside the base locus. Thus the pencil is not connected. The associated weak multinet is a multinet only in case $d=1.$ (Compare with Remark~\ref{refinement}.)

Conversely, if $d=1$ then any $(k,d)$-multinet is a $(k,d)$-net. For every $k$ a $(k,1)$-net is a
pencil of $k$ lines partitioned into singletons. Such an arrangement supports a weak (k/m,m)-multinet for every $m$ dividing $k$. 
From point of view of resonance varieties, $(k,1)$ nets correspond to local
components \cite{Fa97}.
\label{notconnected}
\end{example}

This rank-two example is special.  We conjecture that for rank-three 
pencils the connectedness condition is unnecessary, that is, the
connectedness follows from the other condition of Definition~\ref{Ceva}.

Not every weak multinet arises from a pencil, connected or not. In rank two, consider the multi-arrangement with defining polynomial $x^2y^2(ax-by)(cx+dy).$ The induced weak (3,2)-multinet corresponds to a pencil if and only if $ad-bc=0.$ Here is a rank-three example.

\begin{example}
\label{hessian}
Consider the Hesse pencil $a(x^3+y^3+z^3)-3bxyz$ of Example~\ref{Hesse}.
Let $C_0,C_1,C_2,C_3$ be the four completely reducible fibers, whose components together form the Hessian arrangement \L\ of twelve lines.
The base $\X$ of the pencil consists of the nine inflection points of any of the smooth fibers
of the pencil.
Now define a weak multinet structure on $\L$ with the base locus $\X$
and three classes formed by the irreducible components of respectively
$C_0C_1,C_2^2,C_3^2$. This is a weak $(3,6)$-multinet that is not a multinet.
It is easy to check directly that the pencil $[C_2^2:C_3^2]$ does not contain
$C_0C_1$ as a fiber. It is interesting to note that the reason is that
the four points $$[a:b]=[0:1],[1:1],[1:\xi],[1:\xi^2],$$ 
where $\xi^3=1,$  corresponding to the special fibers of Hesse pencil, cannot be ordered to form a harmonic set.
\end{example}

Our aim in this section is to prove that, by contrast with the previous example,
every multinet arises from a pencil of Ceva type as in Theorem~\ref{pen-net}.  
Specifically, let $(\N,\X)$ 
be a $(k,d)$-multinet on $(\L,m)$ with mutually relatively prime multiplicities,
and let $C_i=\prod_{\ell \in \L_i} \alpha_\ell^{m(\ell)}$ be the product of the defining linear forms (with multiplicities) of lines from the class $\L_i.$ We will show that $C_1,\ldots, 
C_k$ are collinear in the space of degree $d$ curves, i.e., they lie in a pencil, and that this pencil is connected.

First we need to introduce some notation and prove a combinatorial lemma.
Suppose $(\N,\X)$ is a multinet on the multi-arrangement $(\L,m)$. Fix a class $\L_i.$
For every $\ell_0\in\L_i$ and every  $p\in \X$ put
$$\L_i^p=\{\ell\in\L_i \ | \ p \in \ell\},$$
$$\L_i'(\ell_0)=\{\ell\in\L_i \ | \ \ell\cap\ell_0\in\X\},$$
$$\L_i''(\ell_0)=\{\ell\in\L_i \ | \ \ell\cap\ell_0\not\in\X\}.$$ In 
particular, $\ell_0$ is an element of 
$\L_i''(\ell_0)$, but not $\L_i'(\ell_0).$
We will often omit $\ell_0$ from the notation if no ambiguity results.

\begin{lemma}
\label{Q}
Suppose $k: \L \to \Z_{>0}$. Choose $\ell_0\in \L_i$ so that 
$m(\ell_0)/k(\ell_0)\geq m(\ell)/k(\ell)$ for all $\ell\in\L_i-\{\ell_0\}$.
Then we have
\begin{equation}\sum_{\ell\in\L_i}k(\ell)\geq \sum_{p\in\ell_0\cap\X}\sum_{\ell\in\L_i^p}k(\ell).
\label{*}\end{equation}

If $m(\ell_0)/k(\ell_0)>m(\ell)/k(\ell)$ for some $\ell\in\L_i-\{\ell_0\}$ then the inequality (\ref{*}) is strict.
\end{lemma}

\begin{proof}
We apply Remark~\ref{refinement} to build the generalized Cartan matrix $Q$ associated with $\X.$ The submatrix $Q_i$ corresponding to $\L_i$ is indecomposable by Definition~\ref{multinets}(iv).
The vector $u_i=
\sum m(\ell)\omega_{\ell}$ lies in the kernel of $Q_i.$ The row of $Q_i$
corresponding to $\ell_0$ has $s=|\ell_0\cap\X|-1$ on the diagonal,
 -1 in the columns
corresponding to $\ell\in\L_i''\setminus\{\ell_0\}$,
and 0 in all the other columns. This implies
$$sm(\ell_0)=\sum_{\ell\in\L_i'',\ell\not=\ell_0}m(\ell).$$ Multiplying by 
$k(\ell_0)/m(\ell_0)$ and using the assumptions we obtain
\begin{equation}
(s+1)k(\ell_0)=\sum_{\ell\in\L_i''}k(\ell_0)m(\ell)/m(\ell_0)
\leq \sum_{\ell\in\L_i''}k(\ell).
\label{**}
\end{equation}
Now we consider the right-hand side of (\ref{*}). Notice that $\ell_0\in\L_i^p$ for every 
$p\in\ell_0\cap\X$ and the sets $\L_i^p\setminus\{\ell_0\}$ for  $p\in \ell_0\cap \X$  partition 
$\L_i'$. Then, using  (\ref{**})
we obtain the following:
$$\sum_{p\in\ell_0\cap\X}\sum_{\ell\in\L_i^p}k(\ell)=
(s+1)k(\ell_0)+\sum_{\ell\in\L_i'} k(\ell)
\leq \sum_{\ell\in\L_i''}k(\ell)+
\sum_{\ell\in\L_i'}k(\ell)=\sum_{\ell\in\L_i}k(\ell),$$
with strict inequality under the additional condition stated in the lemma.
\end{proof}

\begin{theorem}
\label{pencil}
Suppose $(\N,\X)$ is a multinet on the multi-arrangement $(\L,m).$ Let $C_i=\prod_{\ell \in \L_i}\alpha_\ell^{m(\ell)}.$ Then the pencil of degree $d$
curves determined by any two of $C_1,\ldots, C_k$ contains all of them and is
connected.
\end{theorem}

\begin{proof} Recall that the hypothesis means in particular that the line multiplicities $m(\ell)$ are mutually relatively prime (see Remark~\ref{prime}). Let $\pi: \CP^2 \to \CP^1$ be the pencil generated by $C_1$ and $C_2.$ As a first step we prove there is a fiber $F_i$ of $\pi$  
that is divisible by
$\prod_{\ell \in \L_i} \alpha_\ell,$ for each $i$. By condition (iv) of Definition~\ref{multinets}, we can write $\L_i=\{\ell_1,\ldots, \ell_r\}$ so that, for each $j, \ 1\leq j< r,$
$p_j=\ell_{j+1}\cap \ell_j\not\in \X$, i.e., $p_j$
is not in the base $C_1\cap C_2$ of $\pi$.  
Then there is a unique fiber $F_i$ of $\pi$ passing through $p_1$. 
By Lemma \ref{numerology}(iii), $\ell_1$ contains $d$ points of the base locus 
$C_1\cap C_2$ counting multiplicities. 
Also $\ell_1$ contains $p_1$ by construction. 
The degree $d$ curve $F_i$ contains these $d+1$ points of the line $\ell_1,$ 
counting multiplicities. It follows that  $\ell_1$ is an irreducible
component of $F_i,$ by B\'ezout's theorem. 
Assuming that $\ell_j$
is a component of $F_i$ we have $p_j\in F_i\cap \ell_{j+1},$ whence as above 
$\ell_{j+1}$ is also a component of
$F_i$. Then by induction on $j$ we conclude $\prod_{\ell \in \L_i} \alpha_\ell$ divides $F_i.$ 

Now for each $\ell \in \L_i,$ let $k(\ell)$ be maximal such that $\alpha_\ell^{k(\ell)}$ divides $F_i.$ We have $k(\ell)\geq 1$ for all $\ell \in \L_i$ by the preceding argument.  Let $D_i=\prod_{\ell \in \L_i} \alpha_\ell^{k(\ell)}.$ Then we can write $F_i=D_iE_i$ where $\alpha_\ell$ does not divide $E_i$ for any $\ell \in \L_i.$ 

Next we claim that, for any $\ell_0\in \L_i,$ 
\begin{equation}
\deg(E_i)\geq \sum_{p\in \ell_0\cap \X}\sum_{\ell\in \L_i^p} (m(\ell)- k(\ell)),\label{***}
\end{equation}
 with equality  for all $\ell_0\in \L_i$ if and only if $E_i$ and $D_i$ are disjoint outside the base locus.
Indeed, the degree of $E_i$ can be computed by counting intersections with the line $\ell_0$, with multiplicity. The multiplicity of $E_i$ at any point $p$ of $\X$ is
 $\sum_{\ell\in \L_i^p}(m(\ell) - k(\ell))$ by Lemma~\ref{nodes} applied to $\pi.$ Then we have 
\begin{equation}
\deg(E_i)=\sum_{p\in \ell_0\cap \X} \sum_{\ell\in \L_i^p} (m(\ell)-k(\ell)) +\epsilon,\notag
\end{equation}
where $\epsilon$ is the sum of the multiplicities of $E_i$ at points of $E_i\cap \ell_0$ outside the base locus. This proves the claim.

On the other hand, we have $\deg(E_i)=d-\sum_{\ell \in \L_i}k(\ell),$ and $\sum_{p\in\ell_0\cap \X}\sum_{\ell \in \L_i^p}m(\ell)=\sum_{p\in \ell_0\cap \X} n_p=d,$ by Lemma~\ref{numerology}(iii).
Thus inequality (\ref{***}) is equivalent to 
\begin{equation}
\sum_{\ell \in \L_i} k(\ell) \leq  \sum_{p\in \ell_0\cap \X}\sum_{\ell\in \L_i^p} k(\ell).
\label{****}
\end{equation}
Again, equality holds if and only if $E_i$ and $D_i$ are disjoint outside the base locus.

Combining inequality (\ref{****}) with Lemma~\ref{Q}, we conclude that equality holds in both cases. Then the ratios $m(\ell)/k(\ell)$ are equal to a constant $\rho_i\geq 1$ for all $\ell \in \L_i,$ and  $E_i$ and $D_i$ are disjoint away from the base locus, for every $i.$ We will complete the proof by showing that the first condition implies that $\pi$ is connected. The second condition implies that the proper transforms of $D_i$ and $E_i$ under the blow-up $\varphi: \S \to \CP^2$ of $\CP^2$ at \X\ are disjoint; connectedness of $\pi$ will then imply that $E_i=1$ (i.e., $E_i$ is the ``empty curve"). Then $D_i=F_i$ and $k(\ell)=m(\ell)$ for all $\ell \in \L_i,$ for every $i\geq 3,$ proving the first assertion of the theorem. Furthermore, having $m(\ell)/k(\ell)=1$ for all $\ell\in \L_i,$ for all $i,$ we also conclude that $\pi$ is indeed connected, confirming the second assertion of the theorem.

So, suppose, for every $i\geq 1,$ there is a constant $\rho_i$ such that $m(\ell)/k(\ell)=\rho_i$ for all $\ell \in \L_i.$ Suppose $\pi$ is not connected. By the remarks following Lemma~\ref{nodes}, the rational map $\pi:\CP^2\to \CP^1$ lifts to a regular map $\tilde{\pi}: \S \to\CP^1.$ By the Stein Factorization Theorem (\cite[Cor. III.11.5]{Har77}, see also \cite[p.~556]{GH}), we can write $\tilde{\pi}=f\circ \tilde{\pi}_0$ where $\tilde{\pi}_0: \S \to {\mathcal C}$ is a regular map with connected fibers from \S\ to a curve $\mathcal C$, and $f: {\mathcal C}\to \CP^1$ is a finite regular map. Since $\varphi: \S \to \CP^2$ is birational, $\tilde{\pi}_0$ can be pushed down to a rational map $\pi_0: \CP^2 \to {\mathcal C}$. By assumption, the degree $e$ of $f: {\mathcal C} \to \CP^1$ is greater than 1. The generic fiber of $\pi_0$ is a curve of degree $d'<d,$ nonsingular away from the base locus, and $d=ed'.$

Now, by condition (iv) of Definition~\ref{multinets}, the proper transform $\tilde{D_i}$ of $D_i$ is connected. Since $D_i$ and $E_i$ are disjoint away from the base locus, $\tilde{D_i}$ is in fact a connected component of the proper transform of $F_i=D_iE_i.$ Then $D_i$ is a fiber of $\pi_0.$
Thus the degree $\sum_{\ell \in \L_i} k(\ell)$ of the non-reduced fiber $D_i$ is also equal to $d'.$ (The curve $\mathcal C$ can be embedded in projective space, so that the rational map $\pi_0$ is given by a collection of homogeneous polynomials of degree $d'.$)
By assumption we have $\rho_ik(\ell)=m(\ell)$ for each $\ell \in \L_i.$ Then $ed'=d=\sum_{\ell \in \L_i} m(\ell) =\sum_{\ell \in \L_i} \rho_ik(\ell)=\rho_id',$ 
so $\rho_i=e$ for $i=1,\ldots, k.$ Then $ek(\ell)=m(\ell)$ for every $\ell \in \L,$ contradicting the requirement that the line multiplicities $m(\ell)$ be mutually relatively prime. We conclude that $\pi$ is connected. This completes the proof, by remarks above.
\end{proof}

\begin{corollary} The following are equivalent:
\begin{enumerate}
\item{(i)} \L\ supports a global resonance component of dimension $k-1.$
\item{(ii)} \L\ supports a $(k,d)$-multinet for some $d.$
\item{(iii)} \L\ is the set of components of $k\geq 3$ completely reducible fibers in a Ceva pencil of degree $d$ curves, for some $d.$
\end{enumerate}
\label{bigcor}
\end{corollary}

The preceding result imposes further restrictions on the multiplicities that can appear in a multinet.

\begin{definition}
\label{exponents}
For a completely reducible curve $C=\Pi_{i=1}^q\alpha_i^{m_i},$
the greatest common divisor $e$ of all $m_i$ 
is called {\em the exponent} of $C$. 
\end{definition}

Note $e$ is the exponent of $C$ if and only if $C=D^e$ for some completely reducible curve $D$ of exponent 1. 
If $C_i$ ($i=1,2,\ldots,k$) correspond to the classes $\L_i$ of a multinet then
the exponents $e_1,e_2,\ldots,e_k$ of these curves are mutually relatively prime
according to the agreement in 
Remark \ref{prime}.

\begin{proposition} The exponents of completely reducible fibers of a pencil of Ceva type are pairwise relatively prime.
\label{pairwise}
\end{proposition}

\begin{proof} Suppose $C_1$ and $C_2$ are completely reducible fibers of $\pi$, with exponents $e_1$ and $e_2$. Write $C_i=D_i^{e_i}.$ Suppose $d>1$ is a common divisor of $e_1$ and $e_2$. Then
$$\pi=[C_2:-C_1]=[D_2^{e_2}:-D_1^{e_1}]: \CP^2 \to \CP^1$$
 factors through the mapping $\CP^1\to \CP^1$ given by $[a:b]\mapsto [a^d:b^d].$ Since $d>1$ the latter map has disconnected fibers. Then $\tilde{\pi}=\pi\circ \varphi$ has disconnected fibers, contradicting the definition of Ceva pencil.
\end{proof}

\begin{corollary} The greatest common divisors of the multiplicities within classes of a multinet are pairwise relatively prime.
\end{corollary}

Examples \ref{notconnected} and \ref{hessian} show that these restrictions need not hold for disconnected pencils  or weak multinets.

In fact we have not found any examples of multinets 
with exponents different from 1, and we suspect that none exist. At least we can show that there is only one possible triple of nontrivial exponents.

\begin{theorem}
\label{about exponents}
Suppose $\pi$ is a Ceva pencil with at least three completely reducible fibers
with the exponents greater than 1. Then the only possible values of these
three exponents are 2, 3, and 5, and all other completely reducible fibers have exponent 1. 
\end{theorem}

\begin{proof}
Suppose $C_i$ is a completely reducible fiber of $\pi$, and $C_i=G_i^{e_i}$ with $e_i>1,$ for $i=1,2,3.$
Without loss we can suppose that 
\begin{equation}C_1+C_2=C_3.\label{rel}
\end{equation}
We write 
$$G_i=\prod\limits_{j=1}^{d_i} \alpha_{ij}^{r_{ij}},$$
with the $\alpha_{ij}$ being pairwise nonproportional linear forms in $(x,y,z)$, and
 $d_i$ and $r_{ij}$ being positive integers. We can choose coordinates so that the coefficient of $x$ in every $\alpha_{ij}$ is equal to $1.$ We differentiate both sides of equation~(\ref{rel}) above
with respect to $x$, and conclude that $\alpha_{3,j}^{r_{3,j}-1}$
divides 
\begin{equation}
C_1\sum\limits_{k=1}^{d_1}\frac{e_1r_{1k}}{\alpha_{1k}}+
C_2\sum\limits_{k=1}^{d_2}\frac{e_2r_{2k}}{\alpha_{2k}}\label{mess}\end{equation}
for every $j=1,2,\ldots,d_{3}.$ 
Let $f$ be the numerator of the rational function
$$\sum\limits_{k=1}^{d_1}\frac{e_1r_{1k}}{\alpha_{1k}}-
\sum\limits_{k=1}^{d_2}\frac{e_2r_{2k}}{\alpha_{2k}}$$ when expressed in lowest terms. 
For  $r_{3j}\not =1$, $\alpha_{3,j}^{r_{3,j}-1}$ divides both $C_1+C_2$ and the polynomial ~(\ref{mess}), but does not divide $C_1.$  By eliminating the second term in (\ref{mess}), it follows that $\alpha_{3,j}^{r_{3,j}-1}$ divides $f$ for $j=1,2,\ldots,d_{3}.$ 
Since the degree of $f$ is at most $d_1+d_2-1$ we obtain
$d_1+d_2-1\geq d-d_3,$ where $d$ is the degree of $C_i$, $i=1,2,3$.
Using $d_ie_i\leq d$ it follows that
\begin{equation}
\frac{1}{e_1}+\frac{1}{e_2}+\frac{1}{e_3}>1.\label{halph}
\end{equation}

The latter inequality, together with Proposition~\ref{pairwise}, immediately implies the first statement.  This and Proposition~\ref{pairwise} imply the second statement.\footnote{After we posted a draft of this article,  J. Pereira pointed out to us that the inequality (\ref{halph}) was proved by G.~Halphen around 1884 \cite{Halph} - see also 
 \cite{Neto}.
Besides Pereira managed to prove a similar inequality for $e_2$ and $e_3$ when $e_1=1$.}
\end{proof}
\end{section}

\begin{section}{The Riemann-Hurwitz type formula and fibered line arrangements}

Using the pencil to calculate the euler characteristic of the blowup of $\CP^2$ at the points of $\X$, we obtain a Riemann-Hurwitz type formula which imposes restrictions on parameters of multinets, generalizing analogous results for nets in 
\cite{LY00,Yuz04}. The formula can also be used to show that for some Ceva pencils all singular fibers are completely reducible, thus forming a fibered, $K(\pi,1)$ arrangement. (An arrangement is called $K(\pi,1)$ if its complement is aspherical.)

Let $(\N,\X)$ be a multinet on $(\L,m)$. Recall, for $p\in \X,$ $n_p$ is the number of lines from $\L_i$ containing $p$, independent of $i.$ Let $\pi: \CP^2\to \CP^1$ be the associated pencil, \S\ the blow-up of $\CP^2$ at the points of \X, and $\tilde{\pi}: \S \to \CP^1$ the lift of $\pi$ as in the preceding section. Let $B\subset \CP^1$ be the set of regular values of $\tilde{\pi}$ and $\S_0=\tilde{\pi}^{-1}(B).$

\begin{lemma}\label{bundle}
$\tilde{\pi}: \S_0 \to B$ is the projection of smooth fiber bundle. The fiber is the normalization of a curve of degree $d$ with an ordinary singularity of multiplicity $n_p$ at each $p \in \X$ with $n_p>1.$ 
\end{lemma}

\begin{proof} The first statement is a consequence of the Ehresmann Fibration Theorem \cite{Wall04}. The map $\tilde{\pi}$ is a submersion at every point of $\S_0$ and, as the restriction of a map of compact spaces to a full inverse image, $\tilde{\pi}|_{\S_0}$ is proper. The second statement follows from Lemma~\ref{nodes}(iii), since the generic fiber of $\tilde{\pi}$ is the proper transform of the generic fiber of $\pi.$
\end{proof}

Let $\bar{\X}$ denote the set of intersection points of $\L$ not contained in $\X.$ For $p\in \bar\X,$ let $m_p$ be the multiplicity of $p$ in $\L$. 

\begin{theorem}
Let $(\N,\X)$ be a multinet on $(\L,m)$, and let $\pi: \CP^2 \to \CP^1$ be the associated Ceva pencil. Then 
\begin{equation}
\label{euler}
3+|\X|\geq (2-k)\bigl[3d-d^2+\sum_{p\in\X}(n_p^2-n_p)\bigr]+2|\L|-\sum_{p\in\bar\X}(m_p-1)
\end{equation}
with equality if and only if the classes of $(\N,\X)$ form the only singular fibers of $\pi.$
\label{hurwitz}
\end{theorem}

\begin{proof}
The left-hand side of (\ref{euler}) is precisely the Euler characteristic $\chi(\S)$ since blowing up a variety of dimension 2 at
a point increases its Euler characteristic by 1.
On the other side we have estimated $\chi(\S)$ using the decomposition into fibers of $\tilde{\pi}.$ The euler characteristic of $\S_0$ is determined using the bundle map $\tilde{\pi}: \S_0 \to B$ from Lemma~\ref{bundle}. Suppose $\tilde{\pi}$ has $k'$ singular fibers. Then $\chi(B)=2-k'.$ The euler characteristic of the fiber of $\tilde{\pi}$ is given by the Noether formula $3d-d^2+\sum_{p\in\X} \bigl(n_p^2-n_p\bigr)$ (e.g., see \cite{Wall04}, pp. 157-158). 

There are $k$ special fibers that are products of proper transforms of lines
intersecting only at (not blown up)
points of $\bar \X$. The Euler characteristic of such a
fiber corresponding to a class $\L_i$ is $2|\L_i|-\sum_{p\in\bar\X_i}
(m_p-1)$ where $\bar\X_i=\bar{\X}\cap  \L_i.$ Summing up these contributions, we obtain the right-hand side of (\ref{euler}). If there are no other singular fibers, then $k=k'$ and equality holds. Since the euler characteristic of a singular fiber is strictly greater than the euler characteristic of the generic fiber \cite[pp. 508-510]{GH}, the desired inequality holds in general, and is strict if $k'>k.$
\end{proof}

The inequality (\ref{euler}) is not easy to use. For instance, we
do not know if it implies $k\leq 5$ in general. We can prove this in the following particular case, generalizing the similar statement for nets proved in \cite{LY00}.\footnote{J.~Pereira \cite{Per06} sent us an argument based on his previous work showing that the inequality $k\leq 5$ holds in full generality, without the condition on multiplicities. A proof of this result will appear in his joint paper with the second author.}

\begin{corollary}
Let a line arrangement $\L$ support
a $(k,d)$-multinet $(\N,\X)$
with $d>1$ and $m(\ell)=1$ for every $\ell\in\L$.
Then $k\leq 5$.
\label{five}
\end{corollary}
\begin{proof}
In this case we have $|\L_i|=d$ for each $i$ whence $|\L|=kd$. Since 
any two lines from a class $\L_i$ intersect either at a point from $\X$ or
at a point from $\bar\X_i\subset\bar\X$, we have for every $i$
\begin{equation}
\label{1}
d(d-1)=\sum_{p\in\X}n_p(n_p-1)+\sum_{p\in\bar\X_i}m_p(m_p-1).
\end{equation}
Using $\sum_{p\in\X}n_p^2=d^2$ (Lemma~\ref{numerology}(ii)) this can be transformed to
\begin{equation}
\label{2}
\sum_{p\in\X}n_p=d+\sum_{p\in\bar\X_i}m_p(m_p-1).
\end{equation}
Now using (\ref{2}) we can obtain the following expressions for the left-hand (LHS) and right-hand (RHS) sides of (\ref{euler}):

\begin{equation}
\label{3}
\text{LHS} = 3+|\X| \leq 3+\sum_{p\in\X}n_p=3+d+\sum_{p\in\bar\X_i}m_p(m_p-1)
\end{equation}
for every $i$;
\begin{equation}\label{4}
\text{RHS}=(k-2)(-2d+\sum_{p\in\bar\X_i}m_p(m_p-1))+2kd-\sum_{p\in\bar\X}(m_p-1)
\end{equation}
again for every $i$.
Choosing $i$ that maximizes $\sum_{p\in\bar\X_i}(m_p-1)$ we see that (\ref{euler})
implies 
\begin{equation}\label{5}
\begin{split}
3\geq 3d+(k-3)\sum_{p\in\bar\X_i}m_p(m_p-1)-k\sum_{p\in\bar\X_i}(m_p-1)\\
=3d+(k-3)\sum_{p\in\bar\X_i}(m_p-1)^2-3\sum_{p\in\bar\X_i}(m_p-1)\\
=3d+(k-6)\sum_{p\in\bar\X_i}(m_p-1)^2+3\sum_{p\in\bar\X_i}(m_p-1)(m_p-2)\\
\geq 3d+(k-6)\sum_{p\in\bar\X_i}(m_p-1)^2.
\end{split}
\end{equation}
If $k\geq 6$ then (\ref{5}) cannot hold (recall that $d>1$) which completes the
proof.
\end{proof}

Theorem~\ref{hurwitz} also provides a test for $K(\pi,1)$ arrangements.
From Lemma~\ref{bundle} we obtain the following result.
\begin{corollary} Equality holds in (\ref{euler}) if and only if the restriction of $\pi$ to the complement $M=\CP^2 - (\bigcup \L)$ of $\L$  is a smooth bundle projection with base $B=\CP^1 - (k \ \text{points})$ and fiber a smooth surface with some points removed. In particular, 
\L\ is a $K(\pi,1)$ arrangement.
\end{corollary}
To be precise, the fiber of this bundle map is obtained from a degree $d$ curve with nodes of orders $n_p,$ for $n_p>1$,  by removing $|\X|$ points, including all the nodes.

\begin{definition} A multinet $(\N,\X)$ (or its associated Ceva pencil) is {\em complete} if equality holds in  (\ref{euler}).
\end{definition}

Thus the underlying arrangement of a complete multinet is a $K(\pi,1)$ arrangement.

\begin{example} Let \L\ be the $B_3$ arrangement of Example~\ref{exb3}. Then $d=4, k=3,$ and $|\X|$ consists of seven points with multiplicities $n_p=1,1,1,1,2,2,2.$ Both sides of inequality (\ref{euler}) equal 10, and we conclude that the complement of the $B_3$ arrangement is fibered by quartics with three double points, with the nodes removed (i.e., $\CP^1$'s with six punctures) over $\CP^1$ with three points removed. Therefore the arrangement is $K(\pi,1).$ This arrangement is the complexification of a real simplicial arrangement, and is supersolvable (or fiber-type), each of which implies that the arrangement is $K(G,1)$ - we now have a third proof.

This example generalizes to yield a complete $(3,2r)$-multinet on the reflection arrangement corresponding to the full monomial group $G(r,1,3).$ The corresponding Ceva pencil has three singular fibers: $x^r(y^r-z^r), y^r(z^r-x^r),$ and $z^r(x^r-y^r),$ whose linear factors comprise the arrangement. The coordinate lines $x=0, y=0,$ and $z=0$ are assigned multiplicity $r,$ and all other lines have multiplicity one. There are 3 base points of multiplicity $r$ and $r^2$ base points of multiplicity one. Lines within each class intersect outside the base locus in $r$ double points on the coordinate line, so multinet condition (iv) is satisfied, and the pencil is connected. Both sides of inequality~(\ref{euler}) are equal to $r^2+6,$ so the multinet is complete. It yields a fibering of the complement by curves of degree $2r$ with three nodes of order $r,$ with the nodes removed. These arrangements are supersolvable.
\label{ex2b3}
\end{example}

\begin{example} Equality in (\ref{euler}) also holds for the Ceva arrangement of Example~\ref{Ceva-ex}: $d=k=3, |\X|=9,$ all $n_p=1,$ and there are three triple points outside of $\X,$ the intersections of the lines in each of the three classes.  The complement is fibered by nonsingular curves of degree $d$, over $\CP^2$ with three points removed. Thus the Ceva arrangements are $K(\pi,1)$. They are also supersolvable.

Equality also holds for the Hessian arrangement of Example~\ref{Hesse}: in this case $d=3,k=4, |\X|=9, n_p=1$ for all $p,$ and there are $4\cdot \binom{3}{2}=12$ double points outside of \X. This gives an alternate proof that \L\ is a $K(\pi,1)$ arrangement. The original argument uses the fact that this is the reflection arrangement associated with a Shephard group, and thus the complement is a finite cover of the complement of a complexified simplicial arrangement \cite{OS3}.
\label{Ceva-ex2}
\end{example}

The arrangements above, two infinite families and the Hessian, are the only examples of complete multinets that we have found.

We can use Lemma~\ref{rank2ceva}(iv) as a  ``local test" for the presence of additional singular fibers. Recall $\L_p=\{\ell \in \L \ | \ p\in \ell\}.$

\begin{proposition} Suppose $(\N,\X)$ is a complete multinet. Then, for each $p\in \X,$ $$2n_p-2=\sum_{\ell \in \L_p} (m(\ell)-1).$$ In particular, if $m(\ell)=1$ for all $\ell\in \L,$ then the multinet $(\N,\X)$ is not complete unless it is a net.
\label{localtest}
\end{proposition}

\begin{proof} Suppose the equality fails. Then Lemma~\ref{rank2ceva}(iv) implies there is a fiber of the associated Ceva pencil which has two or more branches at $p$ which are tangent, and doesn't correspond to any class of \N. This fiber remains singular in the blown up surface \S.
\end{proof}

By Proposition~\ref{localtest}, the multinets of Figure~\ref{steiner} is not complete. This is a specialization of the Steiner arrangement of \cite[Figure 3]{Kaw05}. The Steiner arrangement supports a $(3,4)$-net, so Proposition~\ref{localtest} does not apply. But in this case the inequality~\ref{euler} is strict, so the net is not complete. 
The Steiner arrangement is not $K(\pi,1)$, by the ``simple triangle" test of \cite{FaRa86}; we do not know if the arrangement of Figure~\ref{steiner} is $K(\pi,1).$  

\begin{example} A specialization of the Pappus arrangement is pictured in Figure~\ref{pappus}(a). It supports a $(3,3)$-net, with the partition indicated in the figure. Here $d=k=3, |\X|=9,$ and there are six double points and a triple point outside of \X. Lemma~\ref{localtest} yields no conclusion, but Lemma~\ref{hurwitz} detects the existence of additional singular fibers, whose euler numbers sum to two. (The general fiber is a nonsingular cubic, euler number zero.) In this case one can see the hidden fiber in the real picture, as the union of a (vertical) line through three of the base points and a conic through the other six. (This curve has euler number two.)
\end{example}

\begin{example} Here is a family of multinets which has not appeared in the literature. For $r>1$ let $(\L,m)$ be the multi-arrangement with defining polynomial $$Q=\left[(x^r-z^r)(y^r-2^rz^r)\right]\left[(y^r-z^r)(x^r-2^rz^r)\right]\left[(x^r-y^r)z^r\right].$$ Thus all lines have multiplicity one except $z=0$ has multiplicity $r$. In particular the equality of Lemma~\ref{localtest} fails at $p=[1:0:0],$ so the multinet is not complete. The partition indicated by the brackets in the expression of $Q$ defines a $(3,2r)$-multinet on $(\L,m).$ For $r=2$ the arrangement is pictured in Figure~\ref{pappus}(b). The hidden singular fiber does not have a real form, but can be found with some cleverness using {\em Macaulay 2} (thanks to Frank-Olaf Schreyer for this). This arrangement (for $r=2$) is neither simplicial nor supersolvable; nevertheless, somewhat surprisingly in the current context, it is known to be a $K(\pi,1)$ arrangement, by the weight test (\cite{Fa95}, see also \cite{Pa95}). We don't know if this arrangement is $K(\pi,1)$ for $r>2.$
\label{os}
\end{example}

\begin{figure}
\centering
\subfigure[The Pappus arrangement]{\includegraphics[width=5cm]{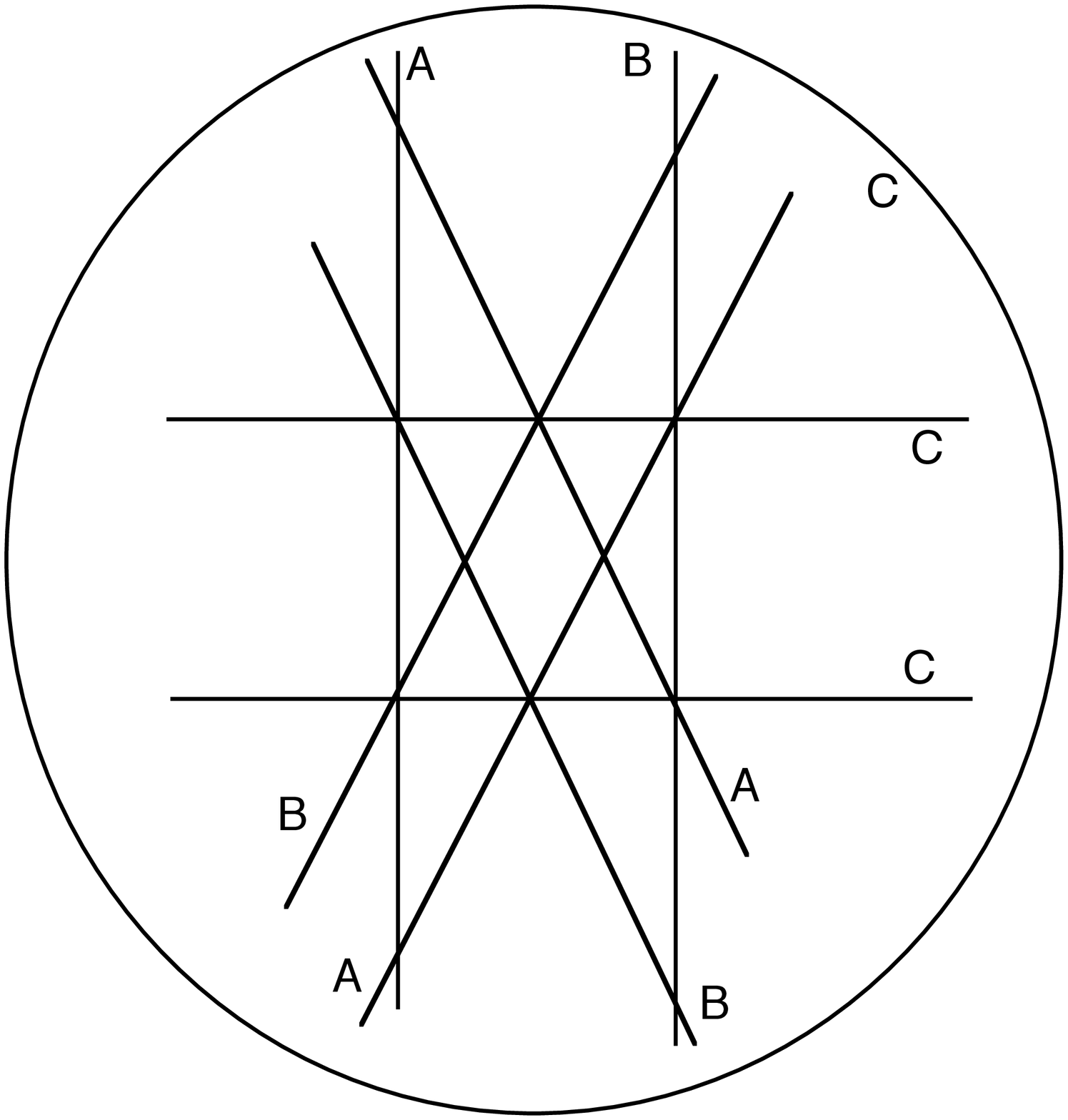}}
\subfigure[A (3,4)-multinet]{\includegraphics[width=5cm]{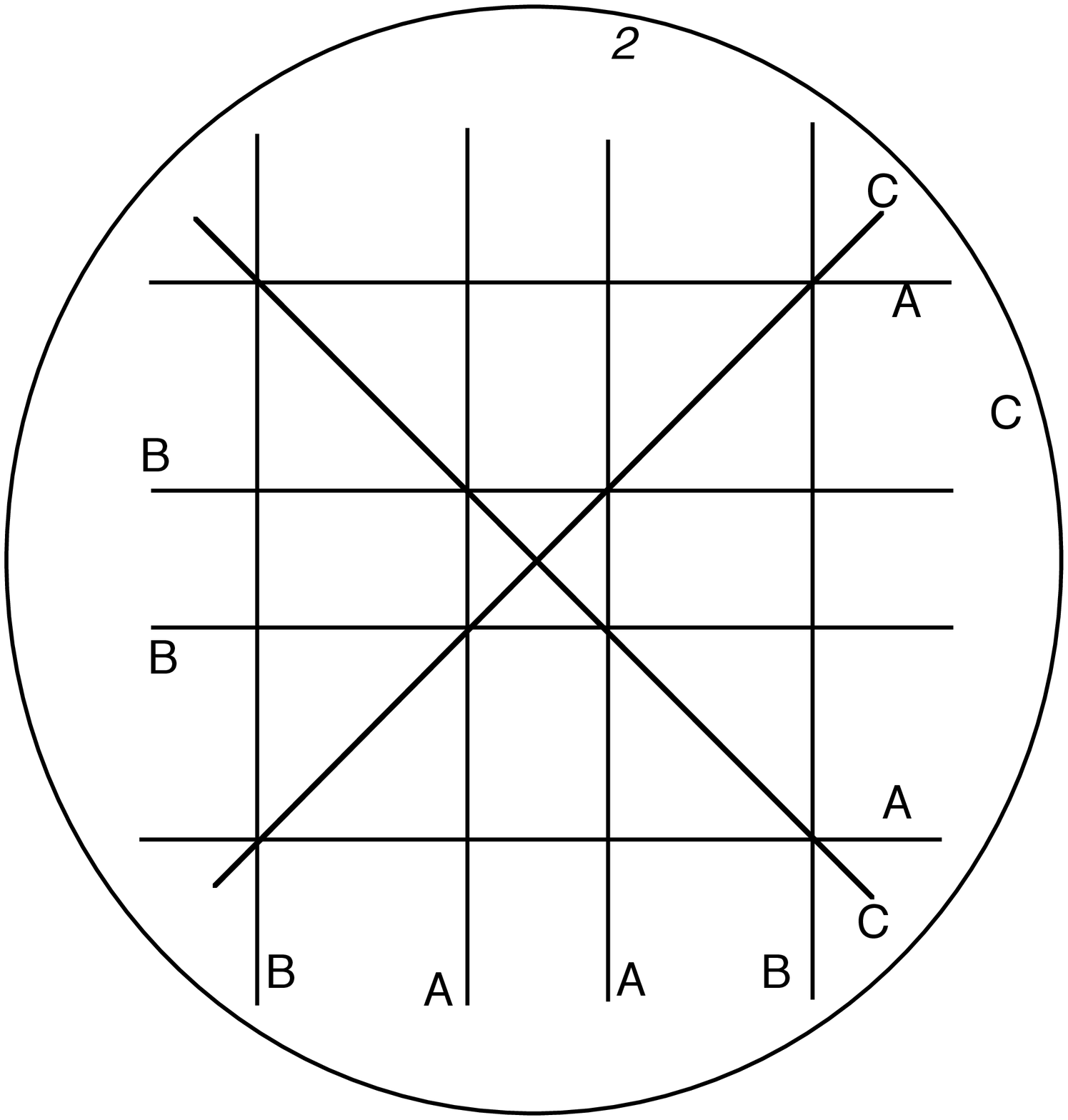}}
\caption{}
\label{pappus}
\end{figure}

\begin{remark} At the level of matroids, every multinet can be obtain from a net by gluing some points and lines (see \cite{Kaw05}). We do not know if this is true on geometric level, i.e. if
every multinet in $CP^2$ can be obtained by a deformation of a net in $CP^2.$
We conjecture that this is true - all examples known to us support this.
If true, this would imply, for instance, that there are no $(k,d)$-multinets with $k\geq 6$
(cf. Corollary~\ref{five}).

For this reason we would like to understand the behavior of Ceva pencils under deformation of the underlying arrangements. For instance, the Pappus configuration degenerates to the $d=3$ Ceva arrangement. The Pappus arrangement is a degeneration of the subarrangement of the Hessian arrangement consisting of three of the classes. The $B_3$ arrangement is a degeneration of the arrangement in Figure~\ref{pappus}(b) (Example~\ref{os}), which is a degeneration of Figure~\ref{steiner}, which in turn is a degeneration of the Steiner arrangement.

To see the effect on the associated pencils of these degenerations, it is convenient to rewrite the right-hand side of inequality (\ref{euler}) using Lemma~\ref{numerology}(ii). We obtain $$3+|\X|\geq (k-2)\left[-3d+\sum_{p\in\X} n_p\right] + 2|\L|-\sum_{p\in\bar\X} (m_p-1).$$ For the sake of argument, let us assume $k=3$; the inequality becomes $$3\geq -3d+\sum_{p\in\X} n_p - |\X| + 2|\L|-\sum_{p\in\bar\X} (m_p-1).$$

In the case of nets, the right-hand side is equal to $3d-\sum_{p\in\bar\X} (m_p-1),$ which is maximized by forcing the most degeneracy outside the base locus. This explains why equality holds for the $d=3$ Ceva arrangement, where all classes are pencils of three lines, but fails for the Pappus arrangement, where two of the classes are in general position. When three classes are in general position, they form three of the four classes of the Hessian, and the pencil contains a fourth singular fiber, which is completely reducible.

In general, keeping the same number of lines, one also wants to bring lines within classes into more special position; this is why the arrangement of Figure~\ref{steiner} comes closer to being fibered than the Steiner arrangement. Bringing lines together, as occurs when Figure~\ref{pappus}(b) degenerates to the $B_3$ arrangement, has the effect of decreasing both sums as well as $|\X|$ and $|\L|,$ the net result of which is difficult to discern.
\end{remark}

In \cite{FaRa86} there appeared examples of (non-supersolvable) fibered $K(\pi,1)$ arrangements closely related to the Ceva arrangements of Examples~\ref{Ceva-ex} and \ref{Ceva-ex2}. The arrangement, called $J_d$ in \cite{FaRa86}, has defining equation $(x^d-y^d)(y^d-z^d)(z^d-x^d)z.$ It is obtained from the Ceva arrangement by adding one line $z=0.$ ($J_2$ is the non-Fano arrangement.) The restriction of the Ceva pencil to the complement of the larger arrangement is also a bundle projection, because the line $z=0$ is transverse to all the regular fibers. (This will be verified below.) Then one can modify the local trivializations of $\pi$ so that they respect the line, and thus restrict to local trivializations of the restricted map.

Here we formulate numerical criteria to determine whether this transversality holds for a given line $\ell_0 \not \in \L$ relative to a fixed Ceva pencil $\pi: \CP^2\to \CP^1$ corresponding to a complete multinet on \L. 

\begin{proposition} Suppose $(\N,\X)$ is a complete $(k,d)$-multinet on $(\L,m)$, and $\ell_0 \not \in \L.$ Then $\ell_0$ is transverse to the regular fibers of $\pi$ if and only if $$2-2d=|(\ell_0\cap (\bigcup \L))-\X|-kd.$$
\label{transverse}
\end{proposition}

\begin{proof} Since $\ell_0$ has multiplicity at most one at each base point, the proper transform $\hat{\ell}_0$ of $\ell_0$ under the blow-up $\varphi: \S \to \CP^2$ is again a line. The restriction of $\tilde{\pi}: \S \to \CP^1$ to $\hat{\ell}_0$ is a $d$-fold branched covering of $\CP^1$ by $\CP^1.$ Branch points arise from tangencies of ${\hat{\ell}}$ with fibers of $\tilde{\pi}.$ The Hurwitz formula (for Riemann surfaces) tells us the number of such tangencies. In this case the formula reads $$2=2d+\sum_{b\in \B} (|\tilde{\pi}^{-1}(b)|-d),$$ where \B\ is the branch locus of $\tilde{\pi}|_{\hat{\ell_0}}$ downstairs. The class $\L_i$ corresponds to a fiber of this map, which has cardinality $t_i=|\ell_0 \cap \bigcup \L_i -\X|,$ because $\hat{\ell_0}$ can only meet the proper transform of $\bigcup \L_i$ outside of the exceptional divisor. Thus, if $t_i<d,$ this fiber contributes $t_i-d$ to the sum on the right.  The line $\ell_0$ is transverse to all the regular fibers of $\pi$ if and only if every point of \B\ corresponds to an $\L_i$ with $t_i<d.$  This in turn is equivalent to $2=2d+\sum_{i=1}^k (t_i-d)=2d+\sum_{i=1}^k t_i - kd= 2d+ |\ell_0\cap (\bigcup \L)-\X|-kd.$
\end{proof}

\begin{example} Let $\ell_0$ be the line $z=0,$ and let $\L$ be the Ceva arrangement, with $\L_1, \L_2$, and $\L_3$ given by the linear factors of $(x^d-y^d), (y^d-z^d),$ and $(z^d-x^d)$ respectively. This defines a complete multinet by Example~\ref{Ceva-ex2}. The line $\ell_0$ meets $\L$ in precisely d+2 points, which are outside the base locus. Since $2-2d=d+2-3d,$ we conclude that $\ell_0$ is transverse to the regular fibers of the associated pencil, which thus induces a fibering of the complement of the arrangement with defining equation $(x^d-z^d)(y^d-z^d)(z^d-x^d).$ The fiber is a curve of degree $d$ with $d^2+d$ points removed, the $d^2$ base points plus $d$ points of $\ell_0.$
\label{jd}
\end{example}

\begin{example} Consider the complete $(3,4)$-multinet on the $B_3$ arrangement in Example~\ref{exb3}. The line $x+y+z=0$ meets the arrangement in six points. Since $2-2\cdot 4=6-3\cdot 4,$ this line is transverse to the regular fibers of the associated Ceva pencil. In fact one can apply Proposition~\ref{transverse} to several lines in succession, so long as they meet inside $\bigcup \L.$ Thus any or all of the lines $x\pm y \pm z=0$ can be added to the arrangement to result in fibered arrangements; all are $K(\pi,1)$ arrangements, fibered by quartics with three double points, with the double points and some additional points removed. 

For the more general $(3,2r)$ multinet of Example~\ref{exb3} one seeks a line which meets the arrangement in $2+2r$ points outside the base locus. The line $x+y+z=0$ does not have that property.
\end{example}

All of the arrangements in the preceding examples are simplicial and/or supersolvable, except the Hessian, which was known to be a $K(\pi,1).$ We have yet to find a ``new" $K(\pi,1)$ arrangement using these methods. On the other hand, all supersolvable arrangements arise through the process described above, starting with the rank-two $(k,1)$-multinet and its associated (complete) Ceva pencil, successively removing lines which miss the base locus and meet inside the arrangement. Here identity of Proposition~\ref{transverse} reads $2-2\cdot 1=d-d\cdot 1,$ that is $0=0.$
\end{section}

\begin{ack}
The authors are grateful to Anatoly Libgober, Rick Miranda, Hiroaki Terao, Dan Cohen, Jose Cogolludo, Hal Schenck, and Frank-Olaf Schreyer, for helpful discussions. We also thank J.~Pereira for informing us of the connection with Halphen's theorem and establishing Corollary~\ref{five} in full generality. Finally we thank the anonymous referee for a careful reading of the manuscript.
\end{ack}

\end{document}